\documentclass[oneside,10pt]{article}          % please do not change
\usepackage[b5paper]{geometry}	    % your paper can be easily printed on a4 or letter paper with enlargenment      
\usepackage{amsfonts,amsmath,latexsym,amssymb,enumerate} % these packages are required
\usepackage{theorem}                % please use one of this two options for theorems
\usepackage{mathrsfs,upref}         % not so essential, but part of journal style
\usepackage{mathptmx}		    % Journal is printed with poscript fonts: 
\usepackage{dea}	            % Journal style, comment only if you find a bug
\newcommand{\mc}{\mathcal{C}}
\newcommand{\mr}{\mathbb{R}}

\newtheorem{theorem}{Theorem}           
\newtheorem{lemma}{Lemma}

\theoremstyle{definition}

\newtheorem{example}{Example}

\newtheorem{remark}{Remark}

\begin{document}

\title[Short title]{On Weakly Nonlinear Boundary Value Problems on Infinite Intervals  }

% Short title is optional, it will appear in running heads.
% It is necessary only if the title is to long to be used in running heads

\author{Benjamin Freedman and Jes\'{u}s Rodr\'{i}guez}

\address{Benjamin Freedman \\ 
Department of Mathematics \\
Box 8205, NCSU, Raleigh, NC 27695-8205\\
USA \\
\email{bnfreedm@ncsu.edu}}

\address{Jes\'{u}s Rodr\'{i}guez \\
Department of Mathematics \\
Box 8205, NCSU, Raleigh, NC 27695-8205 \\
USA \\
\email{rodrigu@ncsu.edu}}

\CorrespondingAuthor{Benjamin Freedman}

%\dedicated{Dedicated to...}                    % Optional

\date{DD.MM.YYYY}                               % Please, write the date of submission


        % AMS-2010 subj class. The list can be found on http://www.ams.org/mathscinet/msc/msc2010.html

%\thanks{This research is supported by ...\par This paper is a lecture that was given at...} 
        % Optional. Only one command thanks is allowed, use \par inside text if you need multiple thanks.       

%\begin{abstract}
%        Abstract of the paper. Please, write few lines of text and try to
%        reduce the number of mathematical formulas.
%\end{abstract}

\begin{abstract}
In this paper, we study weakly nonlinear boundary value problems on infinite intervals. For such problems, we provide criteria for the existence of solutions as well as a qualitative description of the behavior of solutions depending on a parameter. We investigate the relationship between solutions to these weakly nonlinear problems and the solutions to a set of corresponding linear problems.
\end{abstract}

\maketitle

%%%%% END OF TITLE PAGE %%%%%%%%%%%%%%%%%%%%%%%%%%%%%%%%%%%%%%%%%%%%%%

%%%%% BODY OF THE PAPER %%%%%%%%%%%%%%%%%%%%%%%%%%%%%%%%%%%%%%%%%%%%%%
% You should eventually delete (after reading) the rest of the text below %%

\section{Introduction}
\indent The results in this paper pertain to nonlinear boundary value problems on infinite intervals.  We consider problems with weak nonlinearities in both the differential equation and the boundary conditions. We provide a framework which allows us to establish conditions for the existence of solutions and which also enables us to provide a qualitative description of the dependence of solutions on parameters. \\ 
\indent We consider nonlinear boundary value problems on the infinite interval $[0,\infty)$ of the form
       \begin{align}
       x'(t)-A(t)x(t)=h(t)+\varepsilon f(t,x(t)) \label{deq}
       \end{align}
       subject to 
       \begin{align}
       \Gamma(x)=u+\varepsilon \int_0^\infty g(t,x(t)) dt \label{bcdeq}
       \end{align}
       where $A$ is a continuous $n \times n$ matrix-valued function on $[0,\infty)$, $f$ and $g$ are continuously differentiable maps from $\mr^{n+1}$ into $\mr^n$, and $\Gamma$ is a bounded linear map from the space of bounded, continuous functions on $[0,\infty)$ into $\mr^n$. Our main focus will be on the case where the bounded, continuous function $h$ and vector $u \in \mr^n$ are such that the linear problem 
        \begin{align}
       x'(t)-A(t)x(t)=h(t) \label{lin}
       \end{align}
       subject to 
       \begin{align}
       \Gamma(x)=u \label{linbc}
       \end{align}
       has a solution.     \\
\indent In our analysis, we use a scheme somewhat similar to the Lyapunov-Schmidt procedure and results are obtained through an application of the implicit function theorem for Banach spaces. We provide a framework which allows us to determine cases when for $\varepsilon$ sufficiently small in magnitude, \eqref{deq}-\eqref{bcdeq} has solutions which emanate from a particular solution to \eqref{lin}-\eqref{linbc}. \\
\indent There has been extensive literature studying boundary value problems in the context of differential equations on finite intervals. Examples include \cite{freed1}, \cite{ja1}, \cite{ja2}, and \cite{jsuar1}. For results establishing existence of solutions to boundary value problems on infinite intervals the reader is referred to \cite{kartpaper} in the continuous case and  \cite{rodglo}, \cite{rodswe1}, and \cite{rodswe2} in the discrete case. The use of projection methods such as the Lyapunov-Schmidt procedure in the study of boundary value problems is employed in \cite{ces1}, \cite{ces2}, \cite{eth1}, \cite{freed2}, \cite{hale}, \cite{marrod2},\cite{marrod}, \cite{pad},\cite{jfrod}, \cite{roumaw}. and \cite{ura}.

       \section{Differential Equations}
       We use $\mc$ to denote the space of bounded, continuous functions from $[0, \infty)$ into $\mr^n$, and pair this space with the norm $\| x \|_\infty=\sup_{t \geq 0} |x(t)|$. It is clear that $(\mc, \| \cdot \|_\infty)$ is a Banach space. We use  $|\cdot|$ to denote the Euclidean norm on $\mr^n$ and $\| \cdot \|$ for the standard operator norm on the space of $n \times n$ real-valued matrices.  Throughout this section, we assume that $\Gamma: \mc \to \mr^n$ is a bounded linear map and write
       \begin{align}
\| \Gamma \|= \sup_{\|x\|_\infty=1} \left| \Gamma(x) \right|. \nonumber
\end{align}

%and that there exists $\alpha>0$ such that for all $t \geq s \geq 0$,
%       \begin{align}
%      \sup_{t \geq 0} \int_0^t \| \Phi(t) \Phi^{-1}(s)\| < \infty. \nonumber
%       \end{align}

        Let $\Phi(t)$ denote the fundamental matrix for $x'(t)-A(t)x(t)=0$ such that $\Phi(0)=I$ and $\Phi_i$ denote the $i^{th}$ column of $\Phi$ for $1 \leq i \leq n$.  As mentioned in the introduction, our analysis will include a discussion of a set of closely related linear problems. Throughout the paper, the reader will see that conditions we will impose on $A$ guarantee that for any $\psi \in \mc$, $\Phi(\cdot) \int_0^\cdot\Phi^{-1}(s)\psi(s) ds \in \mc$.
        
        We define $\Lambda$ as the $n \times n$ matrix
        \begin{align}
       \Lambda=[\Gamma(\Phi_1(\cdot))| \Gamma(\Phi_2(\cdot))|\cdots | \Gamma(\Phi_n(\cdot))]. \nonumber
       \end{align}
       Note that a function $x \in \mc$ is a solution to
       \begin{align}
       x'(t)-A(t)x(t)=0 \nonumber
       \end{align}
       subject to 
       \begin{align}
      \Gamma(x)=0 \nonumber
       \end{align}
       if and only if $x(0) \in \ker(\Lambda)$. Given $\psi \in \mc$ and $w \in \mr^n$, we know by variation of parameters that any solution to $x'(t)-A(t)x(t)=\psi(t)$ is of the form
       \begin{align}
       x(t)=\Phi(t)x(0)+\Phi(t)\int_0^t \Phi^{-1}(s) \psi(s) ds. \nonumber
       \end{align}
       Imposing the condition that $\Gamma(x)=w$ we get that 
       \begin{align}
       \Lambda x(0)=w-\Gamma\left(\Phi(\cdot) \int_0^\cdot\Phi^{-1}(s)\psi(s) ds  \right). \nonumber
       \end{align}
       Let $p$ denote the dimension of $\ker(\Lambda)$ for some integer $0 \leq p \leq n$. If $p=0$, it is clear that \eqref{lin}-\eqref{linbc} has a unique solution. The bulk of our results concern the case where $p \geq 1$. In this case, we let $W$ be a matrix whose columns form a basis for $[\ker(\Lambda^T)]^\perp$. Note that there exists a solution to the linear boundary value problem
         \begin{align}
       x'(t)-A(t)x(t)=\psi(t) \nonumber
       \end{align}
       subject to 
       \begin{align}
       \Gamma(x)=w \nonumber
       \end{align}
       if and only if 
       \begin{align}
       W^T \left[w -\Gamma\left(\Phi(\cdot) \int_0^\cdot\Phi^{-1}(s)\psi(s) ds  \right)\right]=0. \nonumber
       \end{align}
      Throughout this paper we will mainly be studying the structure of the solution set to \eqref{deq}-\eqref{bcdeq} in the cases when the matrix $\Lambda$ is singular and the corresponding linear problem \eqref{lin}-\eqref{linbc} has a solution, or equivalently where $h$ and $u$ satisfy 
      \begin{align}
       W^T \left[u -\Gamma\left(\Phi(\cdot) \int_0^\cdot\Phi^{-1}(s)h(s) ds  \right)\right]=0. \nonumber
       \end{align}
   
       Based on the discussion above, it is clear that there exists a solution to the nonlinear boundary value problem
       \begin{align}
       x'(t)-A(t)x(t)=h(t)+\varepsilon f(t,x(t)) \nonumber
       \end{align}
       subject to 
       \begin{align}
      \Gamma(x)=u+\varepsilon \int_0^\infty g(t,x(t)) dt \nonumber
       \end{align}
       for $\varepsilon \neq 0$ if there exists $x \in \mc$ and $v \in \ker(\Lambda)$ satisfying
       \begin{align}
       x(t)=\Phi(t)v+ \Phi(t)\int_0^t \Phi^{-1}(s) [h(s)+\varepsilon f(s,x(s)) ]ds \nonumber
       \end{align}
       and 
       \begin{align}
       W^T\left[\int_0^\infty g(t,x(t)) dt-\Gamma\left(\Phi(\cdot) \int_0^\cdot\Phi^{-1}(s)f(s,x(s)) ds \right) \right]=0. \nonumber
       \end{align}
       
       \begin{remark}
It should be observed that the problems we're considering include ones of the form
\begin{align}
\dot{x}(t)-A(t)x(t)= \varepsilon f(t, x(t)) \nonumber
\end{align}
subject to 
\begin{align}
\int_0^\infty B(t) x(t) dt +\sum_{k=0}^\infty C_k x(t_k)=\varepsilon \int_0^\infty g(t, x(t)) dt \nonumber
\end{align}
where $B$ is a function-valued matrix whose entries are integrable functions from $[0,\infty)$ into $\mr^n$. and $C_k$ for $k \geq 0$ is an $n \times n$ matrix with
\begin{align}
\sum_{k=0}^\infty \|C_k\|<\infty. \nonumber
\end{align}. 
\end{remark}

       We now list the following set of conditions which we will impose in our first theorem.
       \begin{enumerate}[I)]
       \item There exists positive constants $K, \alpha$ such that 
       \begin{align}
       \| \Phi(t) \Phi^{-1}(s) \| \leq Ke^{-\alpha(t-s)} \nonumber
       \end{align}
       for all $t \geq s \geq 0$.
       \item For any compact subset $S \subset \mr^n$, $\frac{\partial f}{\partial x}$ is uniformly continuous on $[0,\infty) \times S$ and 
       \begin{align}
       \sup_{t \geq 0}  \left\|  \frac{\partial f}{\partial x}(t,0) \right\|<\infty. \nonumber 
       \end{align}
       
%       \item There exists $x_0 \in \mc^1$ and $\lambda>0$ such that for all $u \in V_\lambda=\{z \in \mc^1: \|z-x_0\|_1<\lambda\}$,
%       \begin{align}
%     \left\|  \frac{\partial f}{\partial x}(t,u(t)) \right\|<m_\lambda. \nonumber
%     \end{align}
     \item For any compact subset $S \subset \mr^n$, $\frac{\partial g}{\partial x}$ is uniformly continuous on $[0,\infty) \times S$ and 
     \begin{align}
     \int_0^\infty \left\| \frac{\partial g}{\partial x}(t,0) \right\| dt <\infty. \nonumber
     \end{align}
     \item For all $h \in \mc$,
     \begin{align}
     \int_0^\infty |g(t,h(t))| dt <\infty. \nonumber
     \end{align}
%     \item For any $\varepsilon>0$, there exists $\delta(\varepsilon)$ such that if $|\zeta_1|<\lambda$ , $|\zeta_2|<\lambda$ and $|\zeta_1-\zeta_2|<\delta$ then
%     \begin{align}
%      \left\| \frac{\partial f}{\partial x}(t,\zeta_1)- \frac{\partial f}{\partial x}(t,\zeta_2) \right\| <\varepsilon\nonumber
%     \end{align}
%     for all $t \geq 0$.
%      \item For any $\varepsilon>0$, there exists $\delta(\varepsilon)$ such that if $|\zeta_1|<\lambda$ , $|\zeta_2|<\lambda$ and $|\zeta_1-\zeta_2|<\delta$ then
%     \begin{align}
%     \left\| \frac{\partial g}{\partial x}(t,\zeta_1)- \frac{\partial g}{\partial x}(t,\zeta_2) \right\| <\varepsilon\nonumber
%     \end{align}
%     for all $t \geq 0$.
     \item There exists an integrable $s: [0,\infty) \to \mr$ satisfying
     \begin{align}
     \left\| \frac{\partial g}{\partial x}(t,x_1)- \frac{\partial g}{\partial x}(t,x_2) \right\| \leq s(t)|x_1-x_2| \nonumber
     \end{align}
     for all $t \geq 0$ and $x_1,x_2 \in \mr^n$.
       \end{enumerate}
       
    Note that for $x \in \mc$, $v \in \ker(\Lambda)$, $\varepsilon \in \mr$, and $t \geq 0$ we have that,
       \begin{align}
      & \left| x(t)-\Phi(t)v-\Phi(t)\int_0^t \Phi^{-1}(s)[h(s)+\varepsilon f(s,x(s))]ds \right| \nonumber\\
       &\leq \|x\|_\infty+\sup_{s \geq 0} \|\Phi(s)\|+\int_0^\infty \|\Phi(t)\Phi^{-1}(s)\||h(s)+\varepsilon f(s,x(s)) | ds \nonumber \\
       &\leq \|x\|_\infty+\sup_{s \geq 0} \|\Phi(s)\|+[\|h\|_\infty+|\varepsilon|\sup_{s \geq 0} |f(s,x(s))|]K\int_0^\infty e^{-\alpha(t-s)} ds \nonumber \\
       &=\|x\|_\infty+\sup_{s\geq 0} \|\Phi(s)\|+[\|h\|_\infty+|\varepsilon|\sup_{s \geq 0} |f(s,x(s))|]K\alpha^{-1}. \nonumber
       \end{align}
      Also observe that 
       \begin{align}
      &\left|W^T\left[\int_0^\infty g(t,x(t)) dt-\Gamma\left(\Phi(\cdot) \int_0^\cdot \Phi^{-1}(s)f(s,x(s)) ds \right) \right]\right| \nonumber \\
      &\leq\|W^T\| \left[\int_0^\infty |g(t,x(t))| dt-\|\Gamma\| \left( \int_0^\infty \|\Phi(t)\Phi^{-1}(s)\| |f(s,x(s))|ds \right)\right] \nonumber \\
      &\leq\|W^T\| \left[\int_0^\infty |g(t,x(t))| dt-\|\Gamma\| \left( \sup_{s \geq 0} |f(s,x(s))| K \int_0^\infty e^{-\alpha(t-s)}ds \right)\right] \nonumber \\
      &=\|W^T\| \left[\int_0^\infty |g(t,x(t))| dt-\|\Gamma\| \left( \sup_{s \geq 0} |f(s,x(s))|K\alpha^{-1}\right)\right]<\infty. \nonumber
       \end{align}
       From this is follows that $H$ given by
       \begin{align}
    H((x,v),\varepsilon)=\left[
  \begin{array}{c}
  H_1((x,v),\varepsilon) \\
  H_2((x,v),\varepsilon) 
  \end{array}
  \right] =\left[
  \begin{array}{c}
  x(t)-\Phi(t)v- \Phi(\cdot)\int_0^\cdot \Phi^{-1}(s) [h(s)+\varepsilon f(s,x(s))] ds \\
    \nonumber \\
W^T\left[\int_0^\infty g(t,x(t)) dt-\Gamma\left(\Phi(\cdot) \int_0^\cdot \Phi^{-1}(s)f(s,x(s)) ds \right) \right]
  \end{array}
  \right] \nonumber
    \end{align}
   is a well-defined map from $\mc \times \ker(\Lambda) \times \mr$ to $\mc \times \mr^p$
       
       Our main result will involve an application of the implicit function theorem for Banach spaces \cite{lang}. This requires continuous Fr\'{e}chet differentiability of $H$. 
       
       In the following lemma, for $i=1,2$ we use $\frac{\partial H_i}{\partial (x,v)}$ to denote the partial (Fr\'{e}chet) derivative of $H_i$ with respect to $(x,v)$.
       \begin{lemma}
       Suppose that $I)-V)$ hold. Then for any $((x,v),\varepsilon) \in \mc \times \ker(\Lambda) \times \mr$, the bounded linear maps $\frac{\partial H_1}{\partial (x,v)}((x,v),\varepsilon)$ and $\frac{\partial H_2}{\partial (x,v)}((x,v),\varepsilon)$ exist and are given by
       \begin{align}
      \left[ \frac{\partial H_1}{\partial (x,v)}((x,v),\varepsilon)\right](\psi,w)(t)=\psi(t)-\Phi(t)w-\varepsilon \left( \Phi(t)\int_0^t \Phi^{-1}(s)\frac{\partial f}{\partial x}(s,x(s))\psi(s)ds \right) \nonumber
      \end{align}
      and
       \begin{align}
      \left[ \frac{\partial H_2}{\partial (x,v)}((x,v),\varepsilon)\right](\psi,w)=W^T\left[\int_0^\infty \frac{\partial g}{\partial x}(t,x(t))\psi(t) dt -\Gamma\left(\Phi(\cdot) \int_0^\cdot\Phi^{-1}(s)\frac{\partial f}{\partial x}(s,x(s))\psi(s) ds \right) \right]. \nonumber
      \end{align}
       Further, $H_1$ and $H_2$  are continuously (Fr\'{e}chet) differentiable.
       \end{lemma}
       \begin{proof}
       
       For $x, \psi \in \mc$ and $v,w \in \ker(\Lambda)$ we have that 
       \begin{align}
       &H_1((x+\psi,v+w),\varepsilon)-H_1((x,v),\varepsilon)-\psi(t)+\Phi(t)w+\varepsilon \left( \Phi(t)\int_0^t \Phi^{-1}(s) \frac{\partial f}{\partial x}(s,x(s))\psi(s) ds \right) \nonumber \\
       &=\varepsilon \left( \Phi(t)\int_0^t \Phi^{-1}(s) \left[f(s,(x+h)(s))-f(s,x(s))-\frac{\partial f}{\partial x}(s,x(s))\psi(s)\right]  ds \right). \nonumber
       \end{align}
       
      For $a,b \in \mr^n$, let $L(a,b)$ denote the straight line segment connecting $a$ and $b$. Note that by the mean value theorem, for all $t \geq 0$ we have that
       \begin{align}
     \left| f(t,(x+\psi)(t))-f(t,x(t)) \right|\leq \sup_{\nu(t) \in L(x(t),(x+\psi)(t))} \left| \frac{\partial f}{\partial x}(t, \nu(t)) \psi(t)\right|  \nonumber  \\
      \text{     and    \hspace{8ex}        }  \left| g(t,(x+\psi)(t))-g(t,x(t)) \right| \leq  \sup_{\zeta(t) \in L(x(t),(x+\psi)(t))} \left| \frac{\partial g}{\partial x}(t, \zeta(t))\psi(t) \right| . \nonumber
       \end{align}
       Then we have that for $t \geq 0$,
        \begin{align}
     & \left|  \left( \int_0^t \Phi(t) \Phi^{-1}(s) \left[f(s,(x+h)(s))-f(s,x(s))-\frac{\partial f}{\partial x}(s,x(s))\psi(s)\right]  ds \right)   \right| \nonumber \\
     &\leq \sup_{\nu(s) \in L(x(s),(x+\psi)(s))} \left\| \left[\frac{\partial f}{\partial x}(s, \nu(s))-\frac{\partial f}{\partial x}(s,x(s))\right] \right\| \left(  \int_0^\infty \left\| \Phi(t) \Phi^{-1}(s) \right\|   ds \right) \|\psi \|_\infty  \nonumber \\
    &\leq  \sup_{\nu(s) \in L(x(s),(x+\psi)(s))} \left\| \left[\frac{\partial f}{\partial x}(s, \nu(s))-\frac{\partial f}{\partial x}(s,x(s))\right] \right\| K\alpha^{-1} \|\psi \|_\infty \nonumber
     \end{align}
     and $\sup_{\nu(s) \in L(x(s),(x+\psi)(s))} \left\| \left[\frac{\partial f}{\partial x}(s, \nu(s))-\frac{\partial f}{\partial x}(s,x(s))\right] \right\| K\alpha^{-1} \rightarrow 0$ as $\|\psi \|_\infty \rightarrow 0$ by $II)$.
 
      We also have that 
      \begin{align}
     &\bigg| H_2((x+\psi,v+w),\varepsilon)-H_2((x,v),\varepsilon) -\nonumber \\
     &W^T \left[\int_0^\infty \frac{\partial g}{\partial x}(t,x(t))\psi(t)dt-\Gamma\left(\Phi(\cdot) \int_0^\cdot\Phi^{-1}(s)\frac{\partial f}{\partial x}(s,x(s)) \psi(s) ds \right) \right] \bigg| \nonumber \\
    &=\bigg|W^T\bigg(\int_0^\infty \left[g(s,(x+\psi)(s))-g(s,x(s))-\frac{\partial g}{\partial x}(s,x(s))\psi(s) \right] ds \nonumber \\
    &- \Gamma\left(\Phi(t) \int_0^t \Phi^{-1}(s)\left[f(s,(x+\psi)(s))-f(s,x(s))-\frac{\partial f}{\partial x}(s,x(s))\psi(s)\right] ds \right) \bigg]  \bigg| \nonumber \\
    &\leq \bigg(\|W^T\| \int_0^\infty \sup_{\zeta(s) \in L(x(s),(x+\psi)(s))} \left\| \frac{\partial g}{\partial x}(s,\zeta(s))-\frac{\partial g}{\partial x}(s,x(s)) \right\| ds \|\psi\|_\infty \nonumber \\
    &+\|W^T \|  \| \Gamma\| \sup_{\nu(s) \in L(x(s),(x+\psi)(s))} \left\| \frac{\partial f}{\partial x}(s,\nu(s))-\frac{\partial f}{\partial x}(s,x(s)) \right\| \int_0^\infty  \left\| \Phi(t) \Phi^{-1}(s) \right\| dt \bigg)\|\psi \|_\infty   \nonumber \\
    &\leq \|W^T\| \bigg( \|s\|_{L^1}\|\psi\|_\infty + \| \Gamma\| \sup_{\nu(s) \in L(x(s),(x+\psi)(s))}  \left\| \frac{\partial f}{\partial x}(s,\nu(s))-\frac{\partial f}{\partial x}(s,x(s)) \right\|K \alpha^{-1}  \bigg) \|\psi \|_\infty   \nonumber
      \end{align}
      where $\| \cdot \|_{L^1}$ denotes the standard norm on $L^1[0,\infty)$. Note that $\|W^T\| \bigg( \|s\|_{L^1} \|\psi\|_\infty
    + \| \Gamma\| \sup_{\nu(s) \in L(x(s),(x+\psi)(s))}\left\| \frac{\partial f}{\partial x}(s,\nu(s))-\frac{\partial f}{\partial x}(s,x(s)) \right\| K \alpha^{-1} \bigg) \rightarrow 0$ as $\|\psi \|_\infty \rightarrow 0$ by $II)$.
      Now we will show that the map
      \begin{align}
      (x,v) \mapsto \frac{\partial H_i}{\partial (x,v)} \nonumber
      \end{align}
      is continuous for $i=1,2$. Note that for $\|\psi \|_\infty=1$,
      
      \begin{align}
      &\left\| \left[ \frac{\partial H_1}{\partial (x,v)}(x_1,v_1) -\frac{\partial H_1}{\partial (x,v)}(x_2,v_2)\right] \psi \right\|_\infty \nonumber \\
      &=\sup_{t \in [0,\infty)} \bigg| \left(\int_0^t \Phi(t)\Phi^{-1}(s) \left[ \frac{\partial f}{\partial x}(s,x_1(s))-\frac{\partial f}{\partial x}(s,x_2(s))\right]\psi(s) ds\right)  \bigg| \nonumber \\
    &\leq  \left\| \frac{\partial f}{\partial x}(s,x_1(s))-\frac{\partial f}{\partial x}(s,x_2(s)) \right\| \left(\int_0^\infty \left\| \Phi(t) \Phi^{-1}(s) \right\| dt\right)  \nonumber \\
    &\leq  K\left\| \frac{\partial f}{\partial x}(s,x_1(s))-\frac{\partial f}{\partial x}(s,x_2(s)) \right\| \alpha^{-1}  \nonumber 
    \end{align} 
    and $ K \left\| \frac{\partial f}{\partial x}(s,x_1(s))-\frac{\partial f}{\partial x}(s,x_2(s)) \right\| \alpha^{-1} \rightarrow 0$ as $\|x_1-x_2\|_\infty \rightarrow 0$. We also have that

      \begin{align}
      &\left| \left[\frac{\partial H_2}{\partial (x,v)}(x_1,v_1) -\frac{\partial H_2}{\partial (x,v)}(x_2,v_2)\right]\psi  \right| \nonumber \\
      &\leq \|W^T\| \bigg( \int_0^\infty   \left\| \frac{\partial g}{\partial x}(s,x_1(s))-\frac{\partial g}{\partial x}(s,x_2(s)) \right\| ds \nonumber \\
    &+ \| \Gamma\| \int_0^\cdot  \| \Phi(\cdot) \Phi^{-1}(s)\| \left\|\frac{\partial f}{\partial x}(s,x_1(s))-\frac{\partial f}{\partial x}(s,x_2(s))\right\|  ds  \bigg] \bigg)  \bigg)  \bigg| \nonumber \\
    &\leq \|W^T\| \bigg(\ \|x_1-x_2\|_\infty \|s\|_{L^1}
    +K\alpha^{-1} \|\Gamma\|  \left\| \frac{\partial f}{\partial x}(s,x_1(s))-\frac{\partial f}{\partial x}(s,x_2(s)) \right\|  \bigg).  \nonumber 
    \end{align}
      Note that $\|W^T\| \bigg(\ \|x_1-x_2\|_\infty  \|s\|_{L^1} 
    + K\alpha^{-1} \|\Gamma\|  \left\| \frac{\partial f}{\partial x}(s,x_1(s))-\frac{\partial f}{\partial x}(s,x_2(s)) \right\|\bigg) \rightarrow 0$ as $\|x_1-x_2\|_\infty \rightarrow 0$, proving our desired result. \qed
       \end{proof}
       
         \begin{remark}
  The most interesting case and the one we will focus mostly on is the case where $\Lambda$ is singular. In this case, solving the nonlinear boundary value problem \eqref{deq}-\eqref{bcdeq} is equivalent to solving the operator equation $H_1((x,v),\varepsilon)=H_2((x,v),\varepsilon)=0$. For the sake of completeness in our analysis it is worth mentioning the case where $\Lambda$ is invertible. If $\Lambda$ is invertible, then \eqref{lin}-\eqref{linbc} has a unique solution and the matrix $W$ does not exist. The nonlinear boundary value problem \eqref{deq}-\eqref{bcdeq} is then equivalent to finding a continuous function $x$ and $v \in \mr^n$ satisfying
  \begin{align}
  x(t)-\Phi(t)v-\Phi(t) \int_0^t \Phi^{-1}(s)[h(s)+\varepsilon f(s,x(s))]ds=0 \nonumber
  \end{align}
where
  \begin{align}
  v=\Lambda^{-1}\left[u+\varepsilon \int_0^\infty g(t, x(t)) dt-\Gamma\left(\Phi(\cdot) \int_0^\cdot \Phi^{-1}(s)[h(s)+\varepsilon f(s,x(s))]ds\right)\right]. \nonumber
  \end{align}
  Define $\Psi: \mc \times \mr^{n+1}\to \mc \times \mr^n$ by $[\Psi_1, \Psi_2]^T$ where
  \begin{align}
  \Psi_1((x,v),\varepsilon)(t)=x(t)-\Phi(t)v-\Phi(t) \int_0^t \Phi^{-1}(s)[h(s)+\varepsilon f(s,x(s))]ds \nonumber
  \end{align}
  and 
   \begin{align}
  \Psi_2((x,v),\varepsilon)(t)=v-\Lambda^{-1}\left[u+\varepsilon \int_0^\infty g(t, x(t)) dt-\Gamma\left(\Phi(\cdot) \int_0^\cdot \Phi^{-1}(s)[h(s)+\varepsilon f(s,x(s))]ds\right)\right]. \nonumber
  \end{align}
  and note that $\Psi((\bar{x},v_0), 0)=0$ where $\bar{x}$ denotes the unique solution to $x'(t)-A(t)x(t)=h(t)$ satisfying $x(0)=v_0$ where 
  \begin{align}
  v_0=\Lambda^{-1}\left[u-\Gamma\left(\Phi(\cdot) \int_0^\cdot \Phi^{-1}(s)h(s)ds \right) \right] . \nonumber
  \end{align}
   Further note that by an analogous argument to the one appearing in the previous lemma, $\Psi$ is continuously differentiable at each point in $\mc \times \mr^{n+1}$ under conditions $I)-V)$ and 
  \begin{align}
  \frac{\partial \Psi}{\partial (x,v)}((\bar{x},v_0),0)[\psi,w]^T=[\psi(\cdot)+\Phi(\cdot)w,w]^T \nonumber
  \end{align}
  which is clearly a bijection from $\mc \times \mr^n$ to $\mc \times \mr^n$. Therefore by the implicit function theorem for Banach spaces, there exists a solution to \eqref{deq}-\eqref{bcdeq} for sufficiently small $\varepsilon$ and those solutions converge uniformly to $\bar{x}$ as $\varepsilon$ goes to $0$.
%   the matrix $W$ doesn't exist and so finding solutions to
%  $x(t)- \Phi(t)\int_0^t \Phi^{-1}(s) [h(s)+\varepsilon f(s,x(s))] ds=0$ is equivalent to establishing solutions to \eqref{deq}-\eqref{bcdeq}. Let $\eta: \mc \times \mr \to \mc$ be defined by 
%  \begin{align}
%  \eta(x,\varepsilon)(t)=x(t)- \Phi(t)\int_0^t \Phi^{-1}(s) [h(s)+\varepsilon f(s,x(s))] ds. \nonumber
%  \end{align}
%  It is clear from the argument appearing in lemma 1 that $\eta$ is continuously Frech\'{e}t differentiable with respect to $x$ for all $(x,\varepsilon) \in \mc \times \mr$.
%  Let $\hat{x}$ be the unique solution to $x'(t)-A(t)x(t)=h(t)$ satisfying $x(0)=0$, it is clear that
%  \begin{align}
%  \eta(\hat{x},0)(t)=\hat{x}(t)- \Phi(t)\int_0^t \Phi^{-1}(s) h(s) ds=0. \nonumber
%  \end{align}
%  for all $t \geq 0$.
%  Also note that the Frech\'{e}t derivative of the map $\eta$ is given by  ,
%  \begin{align}
%  \left[\frac{\partial \eta}{\partial x}(\hat{x},0)\right]=I \nonumber
%  \end{align}
% which is clearly a bijection from $\mc $ onto $\mc$. Therefore, we have by the implicit function for Banach spaces that there exists $\varepsilon_0 > 0$ such that for all $|\varepsilon| < \varepsilon_0$ there exists a unique solution to the nonlinear boundary value problem \eqref{deq}-\eqref{bcdeq} and these solutions emanate from $\hat{x}$. 
         \end{remark}

       Now we shift our focus back to the case where $\Lambda$ is singular. For the sake of notation, for any $y \in \mr^n$ we define the function $x_y(t)=\Phi(t)y+\Phi(t)\int_0^t \Phi^{-1}(s) h(s) ds$. We also write
       \begin{align}
       \frac{\partial H}{\partial (x,v)}=\left[
       \begin{array}{c}
       \frac{\partial H_1}{\partial (x,v)} \nonumber \\
      \nonumber  \\
       \frac{\partial H_2}{\partial (x,v)}
       \end{array}
       \right].
       \end{align}
       
       \begin{theorem} \label{the1}
       Suppose that $I)-V)$ hold and that there exists $y \in \ker(\Lambda)$ such that 
       \begin{align}
       W^T\left[\int_0^\infty g(t,x_y(t))dt -\Gamma\left(\Phi(\cdot) \int_0^\cdot\Phi^{-1}(s)f(s,x_y(s)) ds  \right) \right]=0 \nonumber
       \end{align}
       and $\phi: \ker(\Lambda) \to \mr^p$ given by
        \begin{align}
     \phi(w)=W^T\left[\int_0^\infty \frac{\partial g}{\partial x}(t,x_y(t))\Phi(t) dt-\Gamma\left(\Phi(\cdot) \int_0^\cdot\Phi^{-1}(s)\frac{\partial f}{\partial x}(s,x_y(s))\Phi(s) ds  \right) \right]w  \nonumber
       \end{align}
       is a bijection from $\ker(\Lambda) \subset \mr^n$ onto $\mr^p$. Then there exists $\varepsilon_0$ such that for all $|\varepsilon| \leq \varepsilon_0$, the boundary value problem
        \begin{align}
       x'(t)=A(t)x(t)=h(t)+\varepsilon f(t,x(t)) \nonumber
       \end{align}
       subject to 
       \begin{align}
       \Gamma(x)=u+\varepsilon \int_0^\infty g(t,x(t)) dt. \nonumber
       \end{align}
       is guaranteed a solution $x_\varepsilon$. Moreover $\|x_\varepsilon-x_y\|_\infty \rightarrow 0$ as $\varepsilon \rightarrow 0$.
       \end{theorem}

       \begin{proof}
       We have shown that $H$ is continuously differentiable. Note that $H_1((x_y,y),0)=0=H_2((x_y ,y),0)$. Suppose that $ \frac{\partial H}{\partial (x,v)}((x_y ,y),0)(z, v)=0$. Then $z(t)=\Phi(t)v$ for all $t \geq 0$ and therefore
       \begin{align}
       W^T \left[ \int_0^\infty \frac{\partial g}{\partial x}(s, x_y(s))\Phi(s) ds-\Gamma\left(\Phi(\cdot) \int_0^\cdot\Phi^{-1}(s)\frac{\partial f}{\partial x}(s,x_y(s))\Phi(s) ds  \right) \right] v=0 \nonumber
       \end{align}
       implying that $v=0$. Therefore $\frac{\partial H}{\partial (x,v)}((x_y ,y),0)$ is one-to-one. Let $(\hat{h},\hat{v}) \in \mc \times \mr^p$. Then by assumption there exists a unique $w \in \ker(\Lambda)$ satisfying 
       \begin{align}
       W^T \left[ \int_0^\infty \frac{\partial g}{\partial x}(s, x_y(s))\Phi(s) ds-\Gamma\left(\Phi(\cdot) \int_0^\cdot\Phi^{-1}(s)\frac{\partial f}{\partial x}(s,x_y(s))\Phi(s) ds  \right) \right] w=\hat{v}-v_*. \nonumber
       \end{align}
       where $v_*$ denotes the vector
       \begin{align}
       v_*=W^T \left[ \int_0^\infty \frac{\partial g}{\partial x}(s, x_y(s))\hat{h}(s) ds-\Gamma\left(\Phi(\cdot) \int_0^\cdot\Phi^{-1}(s)\frac{\partial f}{\partial x}(s,x_y(s))\hat{h} ds  \right) \right]. \nonumber
       \end{align}
       Therefore
       \begin{align}
        \left[\frac{\partial H_1}{\partial (x,v)}((x_y, y),0) \right](\hat{h}+\Phi(\cdot)w,w)(t)=\hat{h}(t) \nonumber
        \end{align}
       and 
       \begin{align}
        \left[\frac{\partial H_2}{\partial (x,v)}((x_y, y),0) \right](\hat{h}+\Phi(\cdot)w,w)(t)=(\hat{v}-v_*)+v_*=\hat{v} \nonumber
        \end{align}
       and $ \frac{\partial H}{\partial (x,v)}((x_y ,y),0)$ is a bijection from $\mc \times \ker(\Lambda)$ onto $\mc \times \mr^p$. Our result follows from the implicit function theorem for Banach spaces. \qed
       \end{proof}

       In results up to this point, we assume that $h$ is simply an element of $\mc$. In the following set of results, we investigate problems where we know that $h \in \mc \cap L^1[0,\infty)$. In this case, we impose the following set of conditions.
       \begin{enumerate}[$I'$)]
%       \item There exists $x_0 \in \mc^1$ and $\lambda>0$ such that for all $u \in V_\lambda=\{z \in \mc^1: \|z-x_0\|_1< \lambda\}$,
%       \begin{align}
%       \left\| \frac{\partial f}{\partial x}(t,u(t)) \right\| \leq m_\lambda \nonumber
%       \end{align}
\item There exists positive constant $K$ such that 
       \begin{align}
       \| \Phi(t) \Phi^{-1}(s) \| \leq K \nonumber
       \end{align}
for all $t \geq s \geq 0$.
         \item $\frac{\partial g}{\partial x}$ is uniformly continuous on $[0,\infty) \times \mr^n$ and 
     \begin{align}
     \int_0^\infty \left\| \frac{\partial g}{\partial x}(t,0) \right\| dt <\infty. \nonumber
     \end{align}
     \item For all $h \in \mc$,
     \begin{align}
     \int_0^\infty |g(t,h(t))|dt<\infty. \nonumber
     \end{align}
   \item $\frac{\partial f}{\partial x}$ is uniformly continuous on $[0,\infty) \times \mr^n$ and 
       \begin{align}
       \int_0^\infty  \left\|  \frac{\partial f}{\partial x}(t,0) \right\| dt<\infty. \nonumber 
       \end{align}
%      \item For any $\varepsilon>0$, there exists $\delta(\varepsilon)$ such that if $|\zeta_1|<\lambda$ , $|\zeta_2|<\lambda$ and $|\zeta_1-\zeta_2|<\delta$ then
%     \begin{align}
%     \left\| \frac{\partial g}{\partial x}(t,\zeta_1)- \frac{\partial g}{\partial x}(t,\zeta_2) \right\| <\varepsilon. \nonumber
%     \end{align}
       \item There exists $s \in L^1[0,\infty)$ satisfying
     \begin{align}
     \left\| \frac{\partial g}{\partial x}(t,x_1)- \frac{\partial g}{\partial x}(t,x_2) \right\| \leq s(t)|x_1-x_2| \nonumber
     \end{align}
     for all $t \geq 0$ and $x_1,x_2 \in \mr^n$.
     \item There exists $h_1 \in L^1[0,\infty)$ such that for every compact subset $S$ of $\mr^n$ there exists a constant $C$ satisfying 
     \begin{align}
    | f(t,x)| \leq Ch_1(t) \nonumber
    \end{align}
    for all $t \geq 0$ and $x \in S$ and 
    \begin{align}
    |f(t,x_1)-f(t,x_2)| \leq h_1(t)|x_1-x_2| \nonumber
    \end{align}
    for all $x_1,x_2 \in S$ and $t \geq 0$.
      \item There exists $h_2 \in L^1[0,\infty)$ such that for any compact subset $S \subset \mr^n$,
  \begin{align}
  \left\| \frac{\partial f}{\partial x}(k,x_1)-\frac{\partial f}{\partial x}(k,x_2) \right\| \leq h_2(k)|x_1-x_2| \nonumber
  \end{align}
  for all $t \geq 0$ and $x_1, x_2 \in S$. 
       \end{enumerate}
       Before stating the main theorem in this section, it is worth mentioning for the sake of completeness that if $\Lambda$ is invertible, an analogous argument to the one appearing in remark 2 holds. This is because $\Psi$ is continuously differentiable on $\mc \times \mr^{n+1}$ under conditions $I')-VII')$ and satisfies the conditions of the implicit function theorem at the point $((\bar{x},v_0),0)$ where $\bar{x}$ and $v_0$ are defined the same as in remark 2. Therefore, we can guarantee solutions to $\eqref{deq}-\eqref{bcdeq}$ for $\varepsilon$ sufficiently small and these solutions converge uniformly to $\bar{x}$ as the absolute value of $\varepsilon$ goes to zero.
       \begin{theorem}
       Suppose that $I')-VII')$ hold and that there exists $y \in \ker(\Lambda)$ such that 
       \begin{align}
       W^T\left[\int_0^\infty g(t, x_y(t))dt-\Gamma\left(\Phi(\cdot) \int_0^\cdot\Phi^{-1}(s)f(s,x_y(s)) ds  \right) \right]=0 \nonumber
       \end{align}
       and $\phi: \ker(\Lambda) \to \mr^p$ defined by
        \begin{align}
     \phi(w)=W^T\left[\int_0^\infty \frac{\partial g}{\partial x}(t,x_y(t))\Phi(t) dt-\Gamma\left(\Phi(\cdot) \int_0^\cdot\Phi^{-1}(s)\frac{\partial f}{\partial x}(s,x_y(s))\Phi(s) ds  \right) \right]w \nonumber
       \end{align}
       is a bijection from $\ker(\Lambda) \subset \mr^n$ onto $\mr^p$.
       Then there exists $\varepsilon_0$ such that for all $|\varepsilon| \leq \varepsilon_0$, the boundary value problem
        \begin{align}
       x'(t)-A(t)x(t)=h(t)+\varepsilon f(t,x(t)) \nonumber
       \end{align}
       subject to 
       \begin{align}
       \Gamma(x)=u+\varepsilon \int_0^\infty g(t,x(t)) dt. \nonumber
       \end{align}
       is guaranteed a solution $x_\varepsilon$. Moreover $\|x_\varepsilon-x_y\|_\infty \rightarrow 0$ as $\varepsilon \rightarrow 0$.
       \end{theorem}

       \begin{proof}            We wish to show that $H$ is continuously differentiable under this new set of conditions. Recall that
    \begin{align}
    H_1((x+\psi,v+w),\varepsilon)(t)-H_1((x,v),\varepsilon)(t)-\left[\psi(t)-\Phi(t)w+\varepsilon \left( \Phi(t)\int_0^t \Phi^{-1}(s) \frac{\partial f}{\partial x}(s,x(s))\psi(s) ds \right) \right]  \nonumber \\
    =\varepsilon \left(\Phi(t) \int_0^t \Phi^{-1}(s) \left[ f(s,(x+\psi)(s))-f(s,x(s))-\frac{\partial f}{\partial x}(s,x(s))\psi(s) \right] \right) . \nonumber
    \end{align}
   
    We have that 
    \begin{align}
    &\left\| \Phi(\cdot)\int_0^\cdot \Phi^{-1}(s)\left[ f(s,(x+\psi)(s))-f(s,x(s))-\frac{\partial f}{\partial x}(s,x(s))\psi(s) \right] ds \right\|_\infty \nonumber \\
    &\leq K \int_0^\infty \sup_{\nu(s) \in L(x(s),(x+\psi)(s))}  \left\| \frac{\partial f}{\partial x}(s,\nu(s)) -\frac{\partial f}{\partial x}(s,x(s)) \right\| ds\|\psi \|_\infty  \nonumber \\
    &\leq K \|h_2\|_{L^1} \|\psi \|_\infty^2 \nonumber 
    \end{align}
    and $K \|h_2\|_{L^1} \|\psi\|_\infty \rightarrow 0$ as $\| \psi\|_\infty \rightarrow 0$.
Note also that for $\|\psi \|_\infty=1$,
    \begin{align}
    &\left\| \Phi(\cdot) \left( \int_0^\cdot \Phi^{-1}(s) \left[ \frac{\partial f}{\partial x}(s,x_1(s)) -\frac{\partial f}{\partial x}(s,x_2(s)) \right] \psi(s)ds \right) \right\|_\infty \nonumber \\ 
    &\leq \| \Phi(t) \Phi^{-1}(s) \|  \int_0^\infty \left\| \frac{\partial f}{\partial x}(s,x_1(s)) -\frac{\partial f}{\partial x}(s,x_2(s)) \right\|    ds \nonumber \\
    &\leq K \|h_2\|_{L^1} \|x_1-x_2\|_\infty \rightarrow 0 \nonumber
    \end{align}
    as $\|x_1-x_2\|_\infty \rightarrow 0$.
    We also have that 
      \begin{align}
     &\bigg| H_2((x+\psi,v+w),\varepsilon)-H_2((x,v),\varepsilon) -\nonumber \\
     &W^T \left[\int_0^\infty \frac{\partial g}{\partial x}(t,x(t))\psi(t)dt-\Gamma\left(\Phi(\cdot) \int_0^\cdot\Phi^{-1}(s)\frac{\partial f}{\partial x}(s,x(s)) \psi(s) ds \right) \right] \bigg| \nonumber \\
    &=\bigg|W^T\bigg(\int_0^\infty \left[g(s,(x+\psi)(s))-g(s,x(s))-\frac{\partial g}{\partial x}(s,x(s))\psi(s) \right] ds \nonumber \\
    &- \Gamma\left(\Phi(t) \int_0^t \Phi^{-1}(s)\left[f(s,(x+\psi)(s))-f(s,x(s))-\frac{\partial f}{\partial x}(s,x(s))\psi(s)\right] ds \right) \bigg]  \bigg| \nonumber \\
    &\leq \bigg(\|W^T\| \int_0^\infty \sup_{\zeta(s) \in L(x(s),(x+\psi)(s))} \left\| \frac{\partial g}{\partial x}(s,\zeta(s))-\frac{\partial g}{\partial x}(s,x(s)) \right\| ds \|\psi\|_\infty \nonumber \\
    &+\|W^T \|  \| \Gamma\| \sup_{\nu(s) \in L(x(s),(x+\psi)(s))} \left\| \frac{\partial f}{\partial x}(s,\nu(s))-\frac{\partial f}{\partial x}(s,x(s)) \right\| K \int_0^t \left\| \frac{\partial f}{\partial x}(s,\nu(s))-\frac{\partial f}{\partial x}(s,x(s)) \right\|  dt \bigg)\|\psi \|_\infty   \nonumber \\
    &\leq \|W^T\| \bigg( \|s\|_{L^1}\|\psi\|_\infty + \| \Gamma\| K\|\psi\|_\infty \|h_2\|_{L^1} |  \bigg) \|\psi \|_\infty.   \nonumber
      \end{align}
      and $\|W^T\| \bigg( \|s\|_{L^1}\|\psi\|_\infty + \| \Gamma\| \|\psi\|_\infty \|h_2\|_{L^1} | K  \bigg) \rightarrow 0$ as $\| \psi\|_\infty \rightarrow 0$. Also note that for $\|\psi\|_\infty=1$
      
        \begin{align}
      &\left| \left[\frac{\partial H_2}{\partial (x,v)}(x_1,v_1) -\frac{\partial H_2}{\partial (x,v)}(x_2,v_2)\right]\psi  \right| \nonumber \\
      &\leq \|W^T\| \bigg( \int_0^\infty   \left\| \frac{\partial g}{\partial x}(s,x_1(s))-\frac{\partial g}{\partial x}(s,x_2(s)) \right\| ds \nonumber \\
    &+ \| \Gamma\| \| \Phi(\cdot) \int_0^\cdot\Phi^{-1}(s)\| \left\|\frac{\partial f}{\partial x}(s,x_1(s))-\frac{\partial f}{\partial x}(s,x_2(s))\right\|  ds  \bigg] \bigg)  \bigg)  \bigg| \nonumber \\
    &\leq \|W^T\| \bigg(\ \|x_1-x_2\|_\infty \|s\|_{L^1}
    +\alpha^{-1} \|\Gamma\|  \|x_1-x_2\|_\infty \|h_2\|_{L^1}  \bigg).  \nonumber 
    \end{align}
      It is clear that $\|W^T\| \bigg(\ \|x_1-x_2\|_\infty \|s\|_{L^1}
    +\alpha^{-1} \|\Gamma\|  \|x_1-x_2\|_\infty \|h_2\|_{L^1}  \bigg) \rightarrow 0$ as $\|x_1-x_2\|_\infty \rightarrow 0$.
    Therefore, $H_1$ and $H_2$ is continuously differentiable and so $H$ is as well. It follows that
$H$ satisfies the conditions of the conditions of the implicit function theorem for Banach spaces by an analogous argument to the one appearing in theorem 1. \qed

       \end{proof}
       
       \begin{example}
Consider the boundary value problem
\begin{align}
\dot{x}(t)-Ax(t)= \varepsilon f(t, x(t)) \nonumber
\end{align}
subject to 
\begin{align}
\sum_{k=0}^\infty C_k x(t_k)=\varepsilon \int_0^\infty g(t, x(t)) dt \nonumber
\end{align}
where $x: \mathbb{Z}^+ \to \mr^n$, $f: \mr^3 \to \mr^2$ is twice continuously differentiable, $C_k$ is an $2 \times 2$ real-valued matrix and $t_k \geq 0$ for all $k \geq 0$. We assume that 
\begin{align}
\Lambda=\sum_{k=0}^\infty C_k e^{A t_k} \nonumber
\end{align}
is singular.
Suppose that the matrix $A$ is diagonalizable. That is, there exists an invertible matrix 
\begin{align}
P=\left[
\begin{array}{cc}
p_1 & p_2 \\
p_3 & p_4
\end{array}
\right] \nonumber
\end{align} and diagonal matrix 
\begin{align}
B=\left[
\begin{array}{cc}
\alpha & 0 \\
0 & \beta
\end{array}
\right] \nonumber
\end{align}
satisfying
\begin{align}
A=PBP^{-1}. \nonumber
\end{align}
Therefore, we have that 
\begin{align}
A^k=PB^{k}P^{-1} \nonumber
\end{align}
and so
\begin{align}
e^{At}=P\left[ \sum_{k=0}^\infty \frac{1}{k!} B^k t^k \right]P^{-1}. \nonumber
\end{align}

As mentioned above, we assume that $\Lambda$ is singular, which implies that the second row is a scalar multiple of the first. Suppose that the second row of $\Lambda$ is $\kappa$ times row one for some $\kappa \in \mr$. 
It is clear that $\Lambda$ and $\Lambda^T$ have a one-dimensional kernel and that the kernel of $\Lambda^T$ is spanned by the vector $[-\kappa, 1]^T$. Write $g$ as  $g=[g_1,g_2]$. Suppose that there exists $y \in \ker(\Lambda)$ that satisfies for all $t \geq 0$, 
\begin{align}
0&=f_1(t, e^{At} y)=f_2(t, e^{At} y)=\frac{\partial f_1}{\partial x}(t, e^{At} y)=\frac{\partial f_2}{\partial x}(t, e^{At} y) \nonumber \\
&=g_1(t, e^{At} y)=g_2(t, e^{At} y) \nonumber
\end{align}
and

\begin{align}
-\kappa \int_0^\infty \frac{\partial g_1}{\partial x}(t, e^{At} y) dt \neq \int_0^\infty \frac{\partial g_2}{\partial x}(t, e^{At} y) dt. \nonumber
\end{align}
Under these assumptions, we have
\begin{align}
       W^T\left[\int_0^\infty g(t,e^{tA}y) dt-\sum_{k=0}^\infty C_k e^{A s_k} \int_0^te^{A t_k}f(s,e^{As}y) ds  \right]&=W^T\left[\int_0^\infty (0) dt-\sum_{k=0}^\infty C_k e^{A s_k} \int_0^te^{A t_k}(0)ds  \right]  \nonumber \\
       &=0 \nonumber
       \end{align}
       and that
       \begin{align}
       &\left|W^T\left[\int_0^\infty \frac{\partial g}{\partial x}(t,e^{tA}y)-\sum_{k=0}^\infty C_k e^{At_k}\int_0^t e^{-sA}y\frac{\partial f}{\partial x}(s,e^{sA}y) ds dt \right] \right| \nonumber \\
       &=\left| \int_0^\infty \frac{\partial g_1}{\partial x}(t, e^{At} y) -\kappa \left(\frac{\partial g_2}{\partial x}(t, e^{At} y) \right) dt \right| \nonumber \\
       &\neq 0. \nonumber
       \end{align}
       
%Further suppose that $\gamma$ satisfies 
%\begin{align}
%0=-\beta \left[ \sum_{k=0}^\infty (a_k+kb_k) \sum_{l=0}^{k-1} f_1(l, \gamma, \gamma(l-1))+\sum_{k=0}^\infty b_k \sum_{l=0}^{k-1} [-(l+1)f_1(l,\gamma, \gamma(l-1))+f_2(l,\gamma, \gamma(l-1))] \right] + \nonumber\\ \sum_{k=0}^\infty (c_k+kd_k) \sum_{l=0}^{k-1} f_1(l,\gamma, \gamma(l-1)) +\sum_{k=0}^\infty d_k \sum_{l=0}^{k-1} [-(l+1)f_1(l,\gamma, \gamma(l-1))+f_2(l,\gamma, \gamma(l-1)) ] =\nonumber \\
%=-\beta \left[ \sum_{k=0}^\infty (a_k+kb_k) \sum_{l=0}^{k-1} \frac{ \partial f_1}{\partial x}(l, \gamma, \gamma(l-1))+\sum_{k=0}^\infty b_k \sum_{l=0}^{k-1} [-(l+1)\frac{ \partial f_1}{\partial x}(l,\gamma, \gamma(l-1))+\frac{ \partial f_2}{\partial x}(l,\gamma, \gamma(l-1))] \right] + \nonumber\\ \sum_{k=0}^\infty (c_k+kd_k) \sum_{l=0}^{k-1} \frac{ \partial f_1}{\partial x}(l,\gamma, \gamma(l-1)) +\sum_{k=0}^\infty d_k \sum_{l=0}^{k-1} [-(l+1)\frac{ \partial f_1}{\partial x}(l,\gamma, \gamma(l-1))+\frac{ \partial f_2}{\partial x}(l,\gamma, \gamma(l-1)) ] \nonumber
%\end{align}
Thus for $\varepsilon$ sufficiently small in absolute value, we are guaranteed solutions to the nonlinear boundary value problem above.

%Note
%that for any twice continuously differentiable $f: \mr^3 \to \mr^2$ that $\sum_{k=0}^\infty C_k A^k \sum_{l=0}^{k-1} A^{-(l+1)} f(l, \Phi(l)\gamma)$ can be expressed as: 
%\begin{align}
%\left[
%\begin{array}{c}
%\sum_{k=0}^\infty \hat{a}_k \sum_{l=0}^{k-1} \alpha^{-(l+1)}\left[q_1 f_1(l,\gamma, \gamma(l-1))+q_2 f_2(l,\gamma, \gamma(l-1))\right] \sum_{k=0}^\infty \hat{b}_k \sum_{l=0}^{k-1} [-(l+1)f_1(l,\gamma, \gamma(l-1))+f_2(l,\gamma, \gamma(l-1)) ] \\
%\sum_{k=0}^\infty \hat{c}_k \sum_{l=0}^{k-1} f_1(l,\gamma, \gamma(l-1))+\sum_{k=0}^\infty \hat{d}_k \sum_{l=0}^{k-1} [-(l+1)f_1(l,\gamma, \gamma(l-1))+f_2(l,\gamma, \gamma(l-1)) ]  \nonumber
%\end{array}
%\right].
%\end{align}
%is identical to the second row. \\

%\begin{align}
%W^T \left[ \sum_{k=0}^\infty C_k A^k \sum_{l=0}^{k-1} A^{-(l+1)} f(l, \Phi(l)\gamma) \right] =0 \nonumber
%\end{align}
%Our condition on $g$ implies that 
%\begin{align}
%W^T \sum_{k=0}^\infty g(k, \phi(k) \gamma)&= \sum_{k=0}^\infty \left[ g_1( k, \gamma, \gamma(k-1))- g_2( k, \gamma, \gamma(k-1)) \right]  \nonumber \\
%&=0 \nonumber
%\end{align}

Alternatively, suppose for the problem above that the rows of $\Lambda$ are identical, that $A$ is the matrix
\begin{align}
A=\left[
\begin{array}{cc}
-\frac{1}{2} & 0 \nonumber \\
1 & -\frac{1}{2} \nonumber
\end{array}
\right]
\end{align} and that $f: \mr^3 \to \mr^2$ and $g:\mr^3 \to \mr^2$ are given by
\begin{align}
f(t,x_1,x_2)=
\left[
\begin{array}{c}
\frac{(x_1-e^{-t/2})^2}{t^6} \nonumber \\
\frac{(x_1-e^{-t/2})^2+3(x_2-e^{-t/2}(t+1))^2}{t^8} \nonumber
\end{array}
\right]
\end{align}
and
\begin{align}
g(t,x_1,x_2)=
\left[
\begin{array}{c}
\frac{x_1^2-e^{-t}}{t^2} \nonumber \\
\frac{5(te^{-t/2}-e^{-t/2}-x_2)}{t^2} \nonumber
\end{array}
\right].
\end{align}
Then $y=[1,-1] \in \ker(\Lambda)$ satisfies the conditions imposed in theorem 1. That is, 
\begin{align}
W^T \left[ \int_0^\infty g(t, e^{-t/2}, e^{-t/2}(t-1)) dt+ \sum_{k=0}^\infty C_k e^{At_k} \int_0^t e^{-A(s+1)} f(s, e^{-s/2}, e^{-s/2}(s-1)) ds dt\right] =0 ,\nonumber
\end{align}
and 
\begin{align}
W^T \sum_{k=0}^\infty C_k e^{At_k} \int_0^t e^{-A(s+1)} \frac{\partial f}{\partial x}(s, e^{-s/2}, e^{-s/2}(s-1)) ds dt=W^T \sum_{k=0}^\infty C_k e^{At_k} \int_0^t e^{-A(s+1)} (0) ds dt=0 \nonumber
\end{align}

so we have
\begin{align}
&\left| W^T\left[ \int_0^\infty \frac{ \partial g}{\partial x}(t, e^{-t/2}, e^{-t/2}(t-1)) dt-\sum_{k=0}^\infty C_k e^{At_k} \int_0^t e^{-A(s+1)} \frac{\partial f}{\partial x}(s, e^{-s/2}, e^{-s/2}(s-1)) ds dt \right]\right| \nonumber \\
=&\left| W^T\left[ \int_0^\infty \frac{ \partial g}{\partial x}(t, e^{-t/2}, e^{-t/2}(t-1)) dt \right]\right| \nonumber \\
&=\nonumber \left| W^T \int_0^\infty \left[\frac{ \partial g_1}{\partial x}(t, e^{-t/2}, e^{-t/2}(t-1)) dt - \frac{ \partial g_2}{\partial x}(t, e^{-t/2}, e^{-t/2}(t-1)) \right]dt \right| \nonumber \\
&\neq 0. \nonumber
\end{align}

Therefore, by results in the preceding sections we can guarantee solutions to the nonlinear boundary value problem in this example for $\varepsilon$ sufficiently close to zero.

\end{example}

\end{document}